\documentclass[reqno]{amsart}
\usepackage{amsmath, amssymb, amscd}

\numberwithin{equation}{section}
\theoremstyle{plain}
\newtheorem{thm}{Theorem}

\newtheorem{lem}{Lemma}
\newtheorem{prop}{Proposition}
\newtheorem{Def}{Definition}
\theoremstyle{remark}
\newtheorem{rem}{Remark}

\title{The Askey-Wilson function transform}
\author{Erik Koelink}
\address{Erik Koelink, Technische Universiteit Delft, 
Faculteit Informatietechnologie en Systemen,
Afd. Toegepaste Wiskundige Analyse, Postbus 5031, 2600 GA Delft, The 
Netherlands.}
\email{koelink@twi.tudelft.nl}

\author{Jasper V. Stokman}
\address{Jasper V. Stokman, 
KdV Institute for Mathematics, Universiteit van Amsterdam,
Plantage Muidergracht 24, 1018 TV Amsterdam, The Netherlands.}
\email{jstokman@wins.uva.nl}
\date{April 10th, 2000}
\begin{document}



\begin{abstract}
In this paper we present an explicit (rank one) function transform which
contains several Jacobi-type function transforms and 
Hankel-type transforms as degenerate cases.
The kernel of the transform,
which is given explicitly in terms of basic hypergeometric series, 
thus generalizes the Jacobi function as well as the Bessel
function. The kernel is named the Askey-Wilson function, since it provides
an analytical continuation of the Askey-Wilson polynomial in its
degree.  In this paper we establish the $L^2$-theory of the 
Askey-Wilson function transform, and we explicitely 
determine its inversion formula. 
\end{abstract}

\maketitle


\section{Introduction.}
In \cite{KS2} several explicit function transforms 
are realized as Fourier transforms on the non-compact
quantum $SU(1,1)$ group. 
In this paper we focus on the function theoretic aspects of
the most general Fourier transform on the quantum $SU(1,1)$ group. 
Since the corresponding radial part of the Casimir 
operator is Askey's and Wilson's \cite{AW} second 
order $q$-difference operator $L$, 
we are led to the task of studying its spectral properties
as a symmetric operator with respect to an explicit measure of
unbounded support. This measure naturally arises in \cite{KS2} 
as a Haar weight on the quantum $SU(1,1)$ group.
The kernel of the corresponding generalized Fourier transform
is an eigenfunction of $L$, 
and it provides an analytic continuation
of the Askey-Wilson polynomial in its degree. Consequently, we name
the kernel the Askey-Wilson function.

The Askey-Wilson function satisfies a beautiful symmetry property,
which we refer to as  duality in this paper. It essentially states
that the geometric and the spectral parameter of the Askey-Wilson
function are interchangeable (up to a certain involution on the parameters). 

For a full exploitation of duality 
for the development of the $L^2$-theory of the corresponding 
generalized Fourier transform, 
it is necessary to add another degree of freedom. This can be done 
in the following natural way. 
The Askey-Wilson function can be expanded as a linear
combination of the asymptotically free solutions of $L$ on an
arbitrary $q$-lattice. The corresponding explicit coefficients, 
the so-called $c$-functions, 
can then be  used to define an explicit measure whose support around
infinity lies on the chosen $q$-lattice.
This gives rise to a one parameter family of measures of unbounded
support, in which the extra parameter labels the different 
$q$-lattices. It contains the measure arising 
from the quantum $SU(1,1)$ group as a special case.
 
The generalized Fourier transform with respect to this one parameter
family of measures, and with the Askey-Wilson function as its kernel, is
called the Askey-Wilson function transform. 
In view of its harmonic analytic interpretation 
on the quantum $SU(1,1)$ group,
we may think of the Askey-Wilson function transform as a 
generalization of the classical Jacobi transform, whence in particular
of the Mellin-Fock (or Legendre) transform.
In fact, it was shown in \cite{KS3} that the Askey-Wilson function transform
is on top of an hierarchy consisting of several different
generalizations of the Jacobi and Hankel transform. 
In particular, the $L^2$-theory for the Askey-Wilson function
transform unifies the $L^2$-theory for the spherical Fourier transform
on the (quantum) $SU(1,1)$ group as well as
on the (quantum) group of plane motions.

The main
objective of the paper is to show that the Askey-Wilson function
transform extends to an isometric isomorphism in a natural way,
and that its inverse is the Askey-Wilson function transform
with respect to a different, ``dual'' choice of parameters.

We expect that the present theory on the 
Askey-Wilson function transform 
will turn out to  be the special rank one case 
of a general theory on Macdonald function transforms associated 
with root systems, in analogy with
the polynomial setting (see e.g. \cite{Macd}, \cite{C}).
{}From this point of view, it is natural to expect that
the structure of the Askey-Wilson 
function transform, and eventually also the structure of the general
Macdonald function transform,
will be best understood from an extension of Cherednik's \cite{C}
affine Hecke algebraic approach to the non-polynomial setting.
A first indication of the pivotal role of affine Hecke algebras in the
structure of the Askey-Wilson function transform is the duality, which in the
polynomial setting stems from an involution of the associated double affine
Hecke algebra (see for instance \cite{C2}, and \cite{NS}, \cite{Sahi}
for the rank one setting). We hope that future 
research will shed more light on these matters.

The set-up of the paper is as follows.
In \S 2 we recall the Askey-Wilson polynomials and in \S 3 we define
the Askey-Wilson functions. We emphasize in these sections the
important concept of duality. In \S 4 the asymptotically free
solutions of $L$ are discussed, as well as the corresponding $c$-function
expansion of the Askey-Wilson function. We formulate our main result 
in \S 5, which states that the Askey-Wilson function transform
extends to an isometric isomorphism in a natural way,
and that its inverse is given by the Askey-Wilson function transform
with respect to dual parameters.

The remaining sections occupy the proof of the main result.
Duality plays a pivotal and simplifying role in the proof, 
since it allows us to avoid the machinery of 
spectral analysis of (unbounded) self-adjoint operators (in contrast
with the analysis in some of its degenerate cases, see e.g. \cite{Ka},
\cite{KMU} and \cite{KS2}).
The other main ingredient of the proof 
is the computation of the weak limit of the
Wronskian of the Askey-Wilson functions. 
As an interesting by-product of this computation, we obtain in \S 7 
explicit orthogonality relations for  
``low degree'' Askey-Wilson polynomials with respect to the 
one-parameter family of measures of unbounded support.

{\it Notations and conventions:}
We assume that $0<q<1$ is fixed throughout the paper.
We use the standard notations for basic hypergeometric series,
see for instance \cite{GR}. In particular, we write
$\bigl(a_1,\ldots,a_r;q\bigr)_{c}=\prod_{i=1}^r\bigl(a_i;q\bigr)_c$ with
$\bigl(a;q\bigr)_c=\bigl(a;q\bigr)_{\infty}/\bigl(aq^c;q\bigr)_{\infty}$
and $\bigl(a;q\bigr)_{\infty}=\prod_{i=0}^\infty(1-aq^i)$ 
for (products) of $q$-shifted 
factorials, and we write
\begin{equation*}
{}_{r+1}\phi_r\left(\begin{matrix} a_1, a_2,\ldots, a_{r+1}\\
b_1,b_2,\ldots,b_r \end{matrix};q,z\right)=
\sum_{k=0}^\infty\frac{\bigl(a_1,a_2,\ldots,a_{r+1};q\bigr)_k}{\bigl(q,b_1,
\ldots,b_r;q\bigr)_k}z^k
\end{equation*} 
for the ${}_{r+1}\phi_r$ basic hypergeometric series.
The very well poised ${}_8\phi_7$ is defined by
\begin{equation*}
{}_8W_7(a;b,c,d,e,f;q,z)=
\sum_{k=0}^{\infty}\frac{1-aq^{2k}}{1-a}
\frac{\bigl(a,b,c,d,e,f;q\bigr)_kz^k}{\bigl(q,qa/b,qa/c,qa/d,qa/e,
qa/f;q\bigr)_k}.
\end{equation*}
We use the branch of the square root $\sqrt{\,\cdot\,}$\,
which is positive on ${\mathbb{R}}_{>0}$ throughout the paper.

{\it Acknowledgements:} The second author is 
supported by a fellowship from the Royal
Netherlands Academy of Arts and Sciences (KNAW). Part of the research was done
during the second author's stay at Universit{\'e} Pierre et Marie Curie 
(Paris VI) and Institut de Recherche Math{\'e}matique Avanc{\'e}e (Strasbourg)
in France, supported by a NWO-TALENT
stipendium of the Netherlands Organization for Scientific Research
(NWO) and by the EC TMR network 
``Algebraic Lie Representations'', grant no. ERB FMRX-CT97-0100.


\section{The Askey-Wilson polynomials.}

In this section we briefly recall the basic properties of the
Askey-Wilson polynomials. We formulate these properties in such a way
that the connections with the non-polynomial setting (see section 5)
will be as transparent as possible. 
In particular, we emphasize the important 
concept of duality.

The Askey-Wilson polynomials 
$p_n(x)=p_n(x;a,b,c,d;q)$ ($n\in {\mathbb{Z}}_+$)
are defined by the series expansion
\begin{equation}\label{pn}
p_n(x)={}_4\phi_3\left(\begin{matrix} q^{-n},  q^{n-1}abcd, ax, ax^{-1} \\
 ab, ac, ad \end{matrix}; q,q\right),
\end{equation}
see \cite{AW}.
If we associate dual parameters $\tilde{a}$, $\tilde{b}$, $\tilde{c}$ and 
$\tilde{d}$ with $a,b,c$ and $d$ by the formulas
\begin{equation}\label{dual}
\begin{split}
\tilde{a}&=\sqrt{q^{-1}abcd},\qquad \tilde{b}=ab/\tilde{a}=q\tilde{a}/cd,\\
\tilde{c}&=ac/\tilde{a}=q\tilde{a}/bd,\qquad \tilde{d}=ad/\tilde{a}=
q\tilde{a}/bc,
\end{split}
\end{equation}
then it immediately follows from 
the explicit expression \eqref{pn} that the Askey-Wilson polynomials
satisfy the duality relation
\begin{equation}\label{duality}
p_n(aq^m;a,b,c,d;q)=p_m(\tilde{a}q^n;\tilde{a},\tilde{b},\tilde{c},
\tilde{d};q),\qquad m,n\in {\mathbb{Z}}_+.
\end{equation}
The deeper understanding of duality stems from affine Hecke algebraic 
considerations, see \cite{NS}. 

The Askey-Wilson polynomials $\{p_n\}_{n\in {\mathbb{Z}}_+}$
form a basis of the polynomial algebra ${\mathbb{C}}[x+x^{-1}]$
consisting of eigenfunctions of the Askey-Wilson second order 
$q$-difference operator
\begin{equation}\label{LAW}
\begin{split}
L&=\alpha(x)(T_q-1)+\alpha(x^{-1})(T_q^{-1}-1),\\
\alpha(x)&=\frac{(1-ax)(1-bx)(1-cx)(1-dx)}{(1-x^2)(1-qx^2)},
\end{split}
\end{equation}
where $(T_q^{\pm 1}f)(x)=f(q^{\pm 1}x)$.
The eigenvalue of $L$ corresponding to the Askey-Wilson polynomial
$p_n$ is $\mu(\gamma_n)$, where
$\gamma_n=\tilde{a}q^n$ and 
\begin{equation}\label{eigenvalue}
\mu(\gamma)=-1-\tilde{a}^2+\tilde{a}(\gamma+\gamma^{-1}).
\end{equation}
We remark that 
the three term recurrence relations for the Askey-Wilson polynomials
$p_m(\cdot;\tilde{a},\tilde{b},\tilde{c},\tilde{d};q)$ 
($m\in {\mathbb{Z}}_+$) can be derived from the 
eigenvalue equations $Lp_n=\mu(\gamma_n)p_n$ 
($n\in {\mathbb{Z}}_+$) by applying the duality \eqref{duality}, 
see e.g. \cite{AW} and  \cite{NS}.

For the remainder of this section, 
we restrict our attention to the case that
the parameters $a,b,c$ and $d$ are positive and less than one.
Then Askey and Wilson \cite{AW} proved the orthogonality relations
\begin{equation}\label{orthoAW}
\frac{1}{2\pi i C_0}\int_{x\in {\mathbb{T}}}p_n(x)p_m(x)\Delta(x)\frac{dx}{x}
=\delta_{m,n}\frac{\underset{x=\gamma_0}{\hbox{Res}}\left(
\frac{\widetilde{\Delta}(x)}{x}\right)}
{\underset{x=\gamma_n}{\hbox{Res}}\left(\frac{\widetilde{\Delta}(x)}{x}\right)}
\end{equation}
where $\delta_{m,n}$ is the Kronecker delta and ${\mathbb{T}}$
is the counterclockwise oriented unit circle in the complex plane,
with the weight function given by
\begin{equation}\label{Delta}
\Delta(x)=\frac{\bigl(x^2,1/x^{2};q\bigr)_{\infty}}
{\bigl(ax,a/x,bx,b/x,cx,c/x,dx,d/x;q\bigr)_{\infty}},
\end{equation}
and with $\widetilde{\Delta}(x)$ the weight function $\Delta(x)$
with respect to dual parameters. Here the positive normalization
constant $C_0$ is given by the Askey-Wilson 
integral
\[C_0=\frac{1}{2\pi i}\int_{x\in {\mathbb{T}}}\Delta(x)\frac{dx}{x}=
\frac{2\bigl(abcd;q\bigr)_{\infty}}
{\bigl(q,ab,ac,ad,bc,bd,cd;q\bigr)_{\infty}}.
\]
The orthogonality relations written in the form \eqref{orthoAW} 
exhibit the duality \eqref{duality} of the Askey-Wilson polynomials
on the level of the orthogonality relations, since it expresses 
the quadratic norms explicitly in terms of the dual 
weight function $\widetilde{\Delta}$.
This description of the quadratic norms was proved in 
\cite{NS} using affine Hecke algebraic techniques. 

The above results describe the spectral properties of
the Askey-Wilson second order $q$-difference operator $L$ as a
unbounded symmetric operator on the weighted $L^2$-space with 
respect to the measure $\Delta(x)\frac{dx}{x}$
on the unit circle ${\mathbb{T}}$. For suitable parameter values,
this is directly related to harmonic analysis on the compact
quantum $SU(2)$ group, see e.g. \cite{K}.

In this paper we study the spectral analysis of $L$ as a symmetric operator
with respect to a one-parameter family of measures of unbounded support. 
This setting naturally
arises in the study of harmonic analysis on the non-compact quantum
$SU(1,1)$ group, see \cite{KS2}. Duality, which we have 
emphasized in this section on the polynomial level, 
will play a crucial role in the development
of the corresponding $L^2$-theory.


\section{The Askey-Wilson function.}

In this section we consider a solution of the eigenvalue equation
\begin{equation}\label{eigenvalueequation}
\bigl(Lf\bigr)(x)=\mu(\gamma)f(x),
\end{equation}
which reduces to the Askey-Wilson polynomial for $\gamma=\gamma_n$ 
($n\in {\mathbb{Z}}_+$),
and which enjoys the same symmetry properties as the Askey-Wilson polynomials.

Two linearly independent 
solutions of the eigenvalue equation \eqref{eigenvalueequation}
can be derived from Ismail's and Rahman's 
\cite[(1.11)--(1.16)]{IR} solutions for the
three term recurrence relation of the associated 
Askey-Wilson polynomials. The solutions are given in terms of very well poised
${}_8\phi_7$ series. In this section we consider the solution 
$\phi_{\gamma}(x)=\phi_{\gamma}(x;a;b,c;d;q)$ of \eqref{eigenvalueequation}
given by
\begin{equation}\label{phi}
\begin{split}
\phi_{\gamma}(x)=&\frac{\bigl(qax\gamma/\tilde{d}, qa\gamma/\tilde{d}x;
q\bigr)_{\infty}}
{\bigl(\tilde{a}\tilde{b}\tilde{c}\gamma, q\gamma/\tilde{d},
q\tilde{a}/\tilde{d}, qx/d, q/dx;q\bigr)_{\infty}}\\
&\qquad\times{}_8W_7\bigl(\tilde{a}\tilde{b}\tilde{c}\gamma/q; 
ax, a/x, \tilde{a}\gamma,
\tilde{b}\gamma, \tilde{c}\gamma;q,q/\tilde{d}\gamma\bigr),\qquad
|q/\tilde{d}\gamma|<1.
\end{split}
\end{equation}
This solution of 
\eqref{eigenvalueequation} is also the subject of study in Suslov's papers
\cite{S}, \cite{S2} in which Fourier-Bessel type
orthogonality relations are derived, see also remark \ref{Bessel}.

Applying Bailey's formula \cite[(III.36)]{GR} shows that
\begin{equation}\label{43presentation}
\begin{split}
\phi_{\gamma}(x)=&\frac{1}
{\bigl(bc,qa/d,q/ad;q\bigr)_{\infty}}
{}_4\phi_3\left(\begin{matrix} ax, a/x, \tilde{a}\gamma, \tilde{a}/\gamma\\
 ab, ac, ad \end{matrix}; q,q\right)\\
+&\frac{\bigl(ax, a/x, \tilde{a}\gamma, \tilde{a}/\gamma, qb/d, 
qc/d;q\bigr)_{\infty}}
{\bigl(qx/d, q/dx, q\gamma/\tilde{d}, q/\tilde{d}\gamma, 
ab,ac,bc,qa/d,ad/q;q\bigr)_{\infty}}\\
&\qquad\qquad\times
{}_4\phi_3\left(\begin{matrix} qx/d, q/dx, q\gamma/\tilde{d},
q/\tilde{d}\gamma\\
 qb/d, qc/d, q^2/ad \end{matrix}; q,q\right),
\end{split}
\end{equation}
cf. \cite[(2.8)]{S}, hence
$\phi_{\gamma}(x)$ extends to a meromorphic function in $x$ and $\gamma$
for generic parameters $a,b,c$ and $d$, with possible poles at
$x^{\pm 1}=q^{1+k}/d$ ($k\in {\mathbb{Z}}_+$) and
$\gamma^{\pm 1}=q^{1+k}/\tilde{d}$ ($k\in {\mathbb{Z}}_+$).
It follows from \eqref{43presentation} that
$\phi_{\gamma^{\pm 1}}(x^{\pm 1})
=\phi_{\gamma}(x)$ (all sign combinations possible), and that
$\phi_{\gamma}$ satisfies the duality relation
\begin{equation}\label{dualityfunction}
\phi_{\gamma}(x;a;b,c;d;q)=\phi_x(\gamma;\tilde{a};\tilde{b},
\tilde{c};\tilde{d};q).
\end{equation}
In the remainder of the paper we use the short-hand notation
$\widetilde{\phi}_x(\gamma)$ for the right-hand
side of \eqref{dualityfunction}.

Observe that the meromorphic continuation \eqref{43presentation}
of $\phi_{\gamma}(x)$ implies that 
\begin{equation}\label{redpol}
\phi_{\gamma_n}(x)=\frac{1}
{\bigl(bc,qa/d,q/ad;q\bigr)_{\infty}}
p_n(x),\qquad n\in {\mathbb{Z}}_+,
\end{equation}
since the factor $\bigl(\tilde{a}/\gamma;q\bigr)_{\infty}$
in front of the second ${}_4\phi_3$ vanishes when $\gamma=\gamma_n=
\tilde{a}q^n$ for $n\in {\mathbb{Z}}_+$. In particular, the duality 
\eqref{duality} of the Askey-Wilson polynomials is a special case of
the duality \eqref{dualityfunction} of $\phi_{\gamma}$.

\begin{Def}
The solution $\phi_{\gamma}$ \textup{(}see \eqref{phi} and 
\eqref{43presentation}\textup{)} of the eigenvalue 
equation $(Lf)(x)=\mu(\gamma)f(x)$ is called the Askey-Wilson function.
\end{Def}


\section{The $c$-function expansion.}

We observed in the previous section that the poles of the
Askey-Wilson function $\phi_{\gamma}(x)$ for $|x|\gg 0$ 
lie on the $q$-line $\{dq^{-1-k}\}_{k\in {\mathbb{Z}}_+}$.
In this section we consider the Askey-Wilson function $\phi_{\gamma}(x)$
for $|x|\gg 0$ and $x$ lying on the $q$-line 
$I=\{dtq^k\}_{k\in {\mathbb{Z}}}$, where $t\not\in q^{\mathbb{Z}}$ is an
extra degree of freedom.

Observe that the eigenvalue equation \eqref{eigenvalueequation}
is asymptotically of the form
\begin{equation}\label{asymptequation}
\tilde{a}^2\bigl(f(qx)-f(x)\bigr)+\bigl(f(q^{-1}x)-f(x)\bigr)=\mu(\gamma)f(x)
\end{equation}
when $|x|\rightarrow \infty$.
For generic $\gamma$, 
the asymptotic eigenvalue equation \eqref{asymptequation} has
a basis $\{\Phi_{\gamma}^{free}, \Phi_{\gamma^{-1}}^{free}\}$ of solutions
on the $q$-line $I$, where
\[\Phi_{\gamma}^{free}(dtq^k)=\bigl(\tilde{a}\gamma)^{-k},\qquad
k\in {\mathbb{Z}}.
\]
Furthermore, for generic $\gamma$ there exists a unique solution 
$\Phi_{\gamma}(x)$
of the eigenvalue equation \eqref{eigenvalueequation} on $I$
of the form $\Phi_{\gamma}(x)=\Phi_{\gamma}^{free}(x)g(x)$,
where $g$ has a convergent power series expansion 
around $|x|=\infty$ with constant
coefficient equal to one. $\Phi_{\gamma}$ is the so-called asymptotically
free solution of the eigenvalue equation \eqref{eigenvalueequation}.

An explicit expression for 
$\Phi_{\gamma}$ can be obtained from Ismail's and Rahman's \cite[(1.13)]{IR}
solution of the three term recurrence relation for the associated
Askey-Wilson polynomials. After application of the transformation
formula \cite[(III.23)]{GR} for very well poised ${}_8\phi_7$'s,
it can be expressed as
\begin{equation}
\begin{split}
\Phi_{\gamma}(x)=
&\frac{\bigl(qa\gamma/\tilde{a}x,qb\gamma/\tilde{a}x,qc\gamma/\tilde{a}x,
q\tilde{a}\gamma/dx,d/x;q\bigr)_{\infty}}
{\bigl(q/ax,q/bx,q/cx,q/dx,q^2\gamma^2/dx;q\bigr)_{\infty}}\\
&\times{}_8W_7\bigl(q\gamma^2/dx;q\gamma/\tilde{a},q\gamma/\tilde{d},
\tilde{b}\gamma,\tilde{c}\gamma,q/dx;q,d/x\bigr)\Phi_{\gamma}^{free}(x)
\end{split}
\end{equation}
for $x\in I$ with $|x|\gg 0$.
We now expand the Askey-Wilson function $\phi_{\gamma}(x)$
as a linear combination of the asymptotically free solutions 
$\Phi_{\gamma}(x)$ and $\Phi_{\gamma^{-1}}(x)$ for $x\in I$ with
$|x|\gg 0$. The coefficients in this expansion can be expressed
in terms of the $c$-function $c(\gamma)=c(\gamma;a;b,c;d;q;t)$, which is
defined by
\begin{equation}\label{c}
c(\gamma)=\frac{1}{\bigl(ab,ac,bc,qa/d;q\bigr)_{\infty}\theta(qadt)}
\frac{\bigl(a/\gamma,b/\gamma,
c/\gamma;q\bigr)_{\infty}\theta(\gamma/dt)}
{\bigl(q\gamma/d,1/\gamma^2;q\bigr)_{\infty}}
\end{equation}
where $\theta(x)=(x,q/x;q\bigr)_{\infty}$ is the renormalized
Jacobi theta function.
We call $\widetilde{c}(\gamma)=
c(\gamma;\tilde{a};\tilde{b},\tilde{c};\tilde{d};q;\tilde{t})$
the dual $c$-function, with the dual parameter $\tilde{t}$ defined by
\begin{equation}\label{dualt}
\tilde{t}=1/qadt.
\end{equation}

\begin{prop}\label{cprop}
Let $x\in I$ with $|x|\gg 0$. Then we have the $c$-function expansion
\begin{equation}\label{cexpansion}
\phi_{\gamma}(x)=\widetilde{c}(\gamma)\Phi_{\gamma}(x)+
\widetilde{c}(\gamma^{-1})
\Phi_{\gamma^{-1}}(x)
\end{equation}
for generic $\gamma$.
\end{prop}
\begin{proof}
Apply Bailey's three term recurrence relation \cite[(III.37)]{GR}
with its parameters specialized as follows:
\[
a\rightarrow q\gamma^2/dx,\quad
b\rightarrow q/dx,\quad 
c\rightarrow q\gamma/\tilde{a},\quad 
d\rightarrow q\gamma/\tilde{d},\quad
e\rightarrow \tilde{b}\gamma,\quad 
f\rightarrow \tilde{c}\gamma. 
\]
This gives an expansion of the form \eqref{cexpansion} for explicit
coefficients $\widetilde{c}(\gamma)$ and $\widetilde{c}(\gamma^{-1})$,
which at a first glance still depend on $x\in I$. 
Using the functional 
equation 
\begin{equation}\label{theta}
\theta(q^kx)=q^{-\frac{1}{2}k(k-1)}(-x)^{-k}\theta(x),\qquad k\in {\mathbb{Z}}
\end{equation}
for theta functions, it is easily 
seen that the coefficients are independent of $x$ 
and that they coincide with the dual $c$-functions 
$\widetilde{c}(\gamma^{\pm 1})$ as defined above.
\end{proof}


\section{The Askey-Wilson function transform.}
In this section we define the Askey-Wilson function transform
and we state the main result of this paper. At this stage we need to specify
a particular parameter domain for the five parameters $(a,b,c,d,t)$
in order to ensure positivity of measures.

\begin{Def}
Let $V$ be the set of parameters $(a,b,c,d,t)\in {\mathbb{R}}^{5}$
satisfying the conditions
\begin{equation*}
\begin{split}
&t<0,\qquad\qquad\qquad\qquad 0<b,c\leq a<d/q,\\
&bd,cd\geq q,\qquad\qquad\qquad ab,ac<1.
\end{split}
\end{equation*}
\end{Def}
Observe that $b,c<1$ and $d>q$ for all $(a,b,c,d,t)\in V$.
The domain $V$ is self-dual in the following sense.
\begin{lem}\label{dualdomain}
The assignment 
$(a,b,c,d,t)\mapsto (\tilde{a},\tilde{b},\tilde{c},\tilde{d},
\tilde{t})$
defined by \eqref{dual} and \eqref{dualt}, is an involution on $V$.
\end{lem}
\begin{proof}
Direct verification.
\end{proof}

We fix parameters $(a,b,c,d,t)\in V$ for the remainder of the paper.
For the moment we furthermore 
assume that 
$x\mapsto 1/c(x)c(x^{-1})$ only has simple poles. 
This imposes certain generic conditions on the parameters $(a,b,c,d,t)$, 
which will be removed later by a continuity argument.

It is convenient to renormalize the
function $1/c(x)c(x^{-1})$ as follows,
\begin{equation}\label{Weight}
W(x)=\frac{1}{c(x)c(x^{-1})c_0}
=\frac{\bigl(qx/d,q/dx,x^2,1/x^{2};q\bigr)_{\infty}}
{\bigl(ax,a/x,bx,b/x,cx,c/x;q\bigr)_{\infty}\theta(dtx)\theta(dt/x)},
\end{equation}
where $c_0$ is the positive constant
\begin{equation}\label{c0funct}
c_0=\frac{\bigl(ab,ac,bc,qa/d;q\bigr)_{\infty}^2\theta(adt)^2}{a^2}.
\end{equation}
It follows from \eqref{Weight} that 
\begin{equation}\label{difference}
W(x)=\frac{\theta(dx)\theta(d/x)}{\theta(dtx)\theta(dt/x)}\Delta(x),
\end{equation}
where $\Delta(\cdot)$ is the weight
function \eqref{Delta} for the Askey-Wilson polynomials.
By \eqref{theta}, the quotient
of theta functions in \eqref{difference} is a quasi-constant, i.e. it is
invariant under replacement of $x$ by $qx$. In particular,  
the weight function $W(x)$ differs from $\Delta(x)$ 
only by a quasi-constant factor.

Let $S$ be the discrete subset  
\begin{equation}\label{S}
\begin{split}
S&=\{x\in {\mathbb{C}} \,\, | \,\, 
|x|>1, c(x)=0\}=S_+\cup S_-,\\
S_+&=\{aq^k \,\, | \,\, k\in {\mathbb{Z}}_+,\,\, aq^k>1\},\\
S_-&=\{dtq^k \,\, | \,\, k\in {\mathbb{Z}},\,\, dtq^k<-1\}.
\end{split}
\end{equation} 
We denote $\widetilde{S}$ and $\widetilde{S}_{\pm}$ for the subsets $S$ 
and $S_{\pm}$ with respect to dual parameters.
We define a measure $d\nu(x)$ by
\begin{equation}\label{nu}
\begin{split}
\int f(x)&d\nu(x)=\frac{K}{4\pi i}\int_{x\in {\mathbb{T}}}f(x)
W(x)\frac{dx}{x}\\
&+\frac{K}{2}\sum_{x\in S}f(x)\underset{y=x}{\hbox{Res}}
\left(\frac{W(y)}{y}\right)
-\frac{K}{2}\sum_{x\in S^{-1}}f(x)
\underset{y=x}{\hbox{Res}}
\left(\frac{W(y)}{y}\right),
\end{split}
\end{equation}
where the positive constant
$K$ is given by 
\begin{equation}\label{K}
K=\bigl(ab,ac,bc,qa/d,q;q\bigr)_{\infty}\sqrt{\frac{\theta(qt)\theta(adt)
\theta(bdt)\theta(cdt)}{qabcdt^2}}.
\end{equation}
This particular choice of normalization constant for the measure $\nu$
is justified in theorem \ref{main}.

In view of \eqref{difference}, we can relate the 
discrete masses $\nu\bigl(\{x\}\bigr)(=-\nu\bigl(\{x^{-1}\}\bigr))$ 
for $x\in S_+$ to residues of the weight function $\Delta(\cdot)$,
which were written down explicitly in \cite{AW} (see also \cite[(7.5.22)]{GR}
in order to avoid a small misprint). Explicitly, we obtain 
for $x=aq^k\in S_+$ with $k\in {\mathbb{Z}}_+$ the expression 
\begin{equation}\label{weightplus}
\begin{split}
\nu\bigl(\{aq^k\}\bigr)=&\frac{\bigl(qa/d,q/ad,1/a^2;q\bigr)_{\infty}}
{\bigl(q,ab,b/a,ac,c/a;q\bigr)_{\infty}\theta(adt)\theta(dt/a)}\\
&\quad
\times\frac{\bigl(a^2,ab,ac,ad;q\bigr)_k}{\bigl(q,qa/b,qa/c,qa/d;q\bigr)_k}
\frac{(1-a^2q^{2k})}{(1-a^2)}\frac{K}{2\tilde{a}^{2k}}
\end{split}
\end{equation}
for the corresponding discrete weight. 
For fixed $k\in {\mathbb{Z}}_+$, the right hand side of \eqref{weightplus}
gives the unique continuous extension of
the discrete weight $\nu\bigl(\{aq^k\}\bigr)$ and 
$-\nu\bigl(\{a^{-1}q^{-k}\}\bigr)$ to all parameters 
$(a,b,c,d,t)\in V$ satisfying $aq^k>1$. Furthermore, 
the (continuously extended) 
discrete weight $\nu\bigl(\{aq^k\}\bigr)$ is strictly positive
for these parameter values.

A similar argument can be applied for the discrete weights
$\nu\bigl(\{x\}\bigr)(=-\nu\bigl(\{x^{-1}\}\bigr))$ with $x\in S_-$. 
Explicitly we obtain for $x=dtq^k\in S_-$
with $k\in {\mathbb{Z}}$,
\begin{equation}\label{weighttoinfty}
\begin{split}
\nu\bigl(\{dtq^k\}\bigr)=&\frac{\bigl(qt,q/d^2t;q\bigr)_{\infty}}
{\bigl(q,q,a/dt,b/dt,c/dt,adt,bdt,cdt;q\bigr)_{\infty}}\\
&\times\frac{\bigl(1/t,a/dt,b/dt,c/dt;q\bigr)_{-k}}
{\bigl(q/adt,q/bdt,q/cdt,q/d^2t;q\bigr)_{-k}}
\left(1-\frac{1}{d^2t^2q^{2k}}\right)\frac{K\tilde{a}^{2k}}{2}.
\end{split}
\end{equation}
As for $\nu\bigl(\{x\}\bigr)$ with $x\in S_+$, we use the 
right hand side of \eqref{weighttoinfty} to define 
the strictly positive weight
$\nu\bigl(\{dtq^k\}\bigr)(=-\nu\bigl(\{d^{-1}t^{-1}q^{-k}\}\bigr))$ 
for all $(a,b,c,d,t)\in V$ satisfying  $dtq^k<-1$.

We conclude that the definition of the
measure $\nu$ (see \eqref{nu}) can be extended to arbitrary parameters 
$(a,b,c,d,t)\in V$ using the continuous extensions of its 
discrete weights given above. The resulting measure $\nu$ is a 
positive measure for all $(a,b,c,d,t)\in V$.

\begin{Def}\label{defhilbert}
Let ${\mathcal{H}}$ be the Hilbert space consisting of $L^2$-functions
$f$ with respect to $\nu$ which
satisfy $f(x)=f(x^{-1})$ $\nu$-almost everywhere.
\end{Def}
We write $\widetilde{\nu}$ for the measure $\nu$ with respect to dual 
parameters $(\tilde{a},\tilde{b},\tilde{c},\tilde{d},\tilde{t})$,
and $\widetilde{\mathcal{H}}$ for the associated Hilbert space 
${\mathcal{H}}$. 

Let ${\mathcal{D}}\subset {\mathcal{H}}$
be the dense subspace of functions $f$ with compact support, i.e.
\[ {\mathcal{D}}=\{f\in {\mathcal{H}} \,\, | \,\, f(dtq^{-k})=0,\,\,\,\,
k\gg 0\},
\]
and define
\begin{equation}\label{F}
\bigl({\mathcal{F}}f\bigr)(\gamma)=\int f(x){\overline{\phi_{\gamma}(x)}}
d\nu(x),\qquad f\in {\mathcal{D}}
\end{equation}
for generic $\gamma\in {\mathbb{C}}\setminus \{0\}$. 

\begin{rem}
Observe that the analytic continuation \eqref{43presentation} of 
$\phi_{\gamma}(x)$ is not defined for parameters $(a,b,c,d,t)\in V$
satisfying  $\theta(ad)=0$. These apparent poles can be removed
in view of the original definition \eqref{phi} for 
$\phi_{\gamma}(x)=\phi_{\gamma^{-1}}(x)$
in terms of very well poised ${}_8\phi_7$'s
(observe that $q/\tilde{d}<1$, so that
either $|q\gamma/\tilde{d}|<1$ or $|q/\tilde{d}\gamma|<1$ for 
$\gamma\in {\mathbb{C}}\setminus \{0\}$, hence $\phi_{\gamma}$ can be
expressed in terms of the original definition \eqref{phi}
for generic $\gamma$). In particular, the transform ${\mathcal{F}}$
is well defined for parameters satisfying $\theta(ad)=0$.
\end{rem}
We write ${\widetilde{\mathcal{D}}}\subset \widetilde{{\mathcal{H}}}$ 
(respectively $\widetilde{\mathcal{F}}$) for the dense subspace ${\mathcal{D}}$
(respectively the function transform $\mathcal{F}$)
with respect to dual parameters 
$(\tilde{a},\tilde{b},\tilde{c},\tilde{d},\tilde{t})$.
The main theorem of this paper can now be stated as follows.
\begin{thm}\label{main}
Let $(a,b,c,d,t)\in V$.
The transform ${\mathcal{F}}$ extends to an isometric isomorphism
${\mathcal{F}}: {\mathcal{H}}\rightarrow \widetilde{{\mathcal{H}}}$
by continuity. The inverse of ${\mathcal{F}}$ is given by 
${\widetilde{\mathcal{F}}}$.
\end{thm}

We will discuss the proof of theorem \ref{main}
in detail in the remaining sections.

\begin{Def}
The isometric isomorphism ${\mathcal{F}}: {\mathcal{H}}\rightarrow
\widetilde{\mathcal{H}}$ is called the Askey-Wilson function transform.
\end{Def}

Theorem \ref{main} gives a simultaneous generalization of the
$L^2$-theory for the Mellin-Fock transform and for
the (rank one) Hankel transform, 
see \cite{KS2} and \cite{KS3}.
In fact, the Askey-Wilson function transform is on top of
an hierarchy of several different generalizations of the Mellin-Fock 
and Hankel transform, in analogy with the polynomial setting.
This point of view is emphasized in \cite{KS3}.  
In particular, the Askey-Wilson function generalizes the Jacobi function
as well as the Bessel function to the level of very well poised 
${}_8\phi_7$ series.

Theorem \ref{main} may be partially reformulated in terms
of orthogonality relations for ``low degree'' 
Askey-Wilson polynomials with respect to the measures
$\nu$ ($t<0$). 
We say more about
this point of view in section 7, where we prove the orthogonality relations
for the Askey-Wilson functions $\phi_{\gamma}$ ($\gamma\in \widetilde{S}$)
with respect to the measure $\nu$ as an intermediate step of the proof
of theorem \ref{main}.


\section{The Wronskian.}\label{WWW}

We let $\chi_k\in {\mathcal{D}}$ for $k\in {\mathbb{Z}}$ with $k\ll 0$
be the characteristic function which is zero at $dtq^m, d^{-1}t^{-1}q^{-m}$
for $m\in {\mathbb{Z}}$ with $m<k$ and which is equal to one otherwise.
The sequence $\{\chi_k\}_{k}$ of characteristic functions 
is an approximation of the unit, 
in the sense that $\chi_kf\rightarrow f$ 
in ${\mathcal{H}}$ as $k\rightarrow -\infty$ for all $f\in {\mathcal{H}}$.

The starting point for the proof of theorem \ref{main} is the following
simple observation.
\begin{lem}\label{rewrite}
Let $f,g\in \widetilde{\mathcal{D}}$ and $k\ll 0$, then
\begin{equation*}
\begin{split}
\int\chi_k(x)\bigl(\widetilde{\mathcal{F}}f\bigr)(x)
&\overline{\bigl(\widetilde{\mathcal{F}}g\bigr)(x)}d\nu(x)=\\
=&\iint f(\gamma)\overline{g(\gamma^\prime)}
\left(\int \chi_k(x)\phi_{\gamma}(x)\phi_{\gamma^\prime}(x)
d\nu(x)\right)
d\bigl(\widetilde{\nu}\times\widetilde{\nu}\bigr)(\gamma,\gamma^\prime).
\end{split}
\end{equation*}
\end{lem}
\begin{proof}
We substitute the definition of the dual Askey-Wilson function transform
$\widetilde{\mathcal{F}}$ into the left hand side of the desired
identity, and we use that
\[\overline{\widetilde{\phi}_x(\gamma)}=\widetilde{\phi}_x(\gamma)=
\phi_{\gamma}(x)
\]
for $x\in \hbox{supp}(\nu)$ and $\gamma\in \hbox{supp}(\widetilde{\nu})$,
see \eqref{dualityfunction} for the second equality.
We arrive at
\begin{equation*}
\begin{split}
\int\chi_k(x)\bigl(\widetilde{\mathcal{F}}f\bigr)(x)
&\overline{\bigl(\widetilde{\mathcal{F}}g\bigr)(x)}d\nu(x)=\\
=&\int \chi_k(x) \left(\int f(\gamma)\phi_{\gamma}(x)d\widetilde{\nu}(\gamma)
\right)\left(\int \overline{g(\gamma^\prime)}\phi_{\gamma^\prime}(x)
d\widetilde{\nu}(\gamma^\prime)\right)d\nu(x).
\end{split}
\end{equation*}
Since $f,g\in \widetilde{\mathcal{D}}$, all integrations are over
compact subsets, so we may use Fubini's theorem to interchange the order
of integration. This yields the desired result.
\end{proof}

The proof of theorem \ref{main} now hinges on the explicit evaluation
of the weak limit as $k\rightarrow -\infty$ of the function
\begin{equation}\label{kernel}
(\gamma,\gamma^\prime)\mapsto 
\int\chi_k(x)\phi_{\gamma}(x)\phi_{\gamma^\prime}(x)
d\nu(x)
\end{equation}
with respect to the product measure $\widetilde{\nu}\times \widetilde{\nu}$.

As a first step, we express the integral
\eqref{kernel} in terms of the Wronskian of the Askey-Wilson functions
$\phi_{\gamma}$ and $\phi_{\gamma^\prime}$. Here the Wronskian
$\lbrack f,g\rbrack(x)$ for two functions
$f, g: I \rightarrow {\mathbb{C}}$ 
is defined by
\begin{equation}\label{wronskian}
\lbrack f,g\rbrack(x)=
2\,\nu\bigl(\{x\}\bigr)\alpha(x)\bigl(f(x)g(qx)-f(qx)g(x)\bigr),\qquad
x\in S_-,
\end{equation}
where $\alpha(x)$ is given by \eqref{LAW}.

\begin{prop}\label{propwronskian}
Let $(a,b,c,d,t)\in V$. For generic spectral parameters
$\gamma,\gamma^\prime\in {\mathbb{C}}$ 
such that $\mu(\gamma)\not=\mu(\gamma^\prime)$, we have
\[ \int\chi_k(x)\phi_{\gamma}(x)\phi_{\gamma^\prime}(x)d\nu(x)=
\frac{\lbrack \phi_{\gamma}, \phi_{\gamma^\prime}\rbrack(dtq^{k-1})}
{\mu(\gamma)-\mu(\gamma^\prime)}
\]
for $k\in {\mathbb{Z}}$ sufficiently negative.
\end{prop}
\begin{proof}
The proof of the desired identity simplifies when we slightly perturb
the parameters $(a,b,c,d,t)$. The proof for parameters
$(a,b,c,d,t)\in V$ can then be derived from the perturbed case
using the fact that the left and right hand side of
the identity depend continuously on 
the parameters $(a,b,c,d,t)$.

Let us indicate one class of possible perturbations of 
the parameters which is sufficient for our purposes.
Let $(\alpha,\beta,\gamma,\delta,t)\in V$ with 
$\alpha,\delta,\alpha\delta, -\delta t\not\in q^{\mathbb{Z}}$.
Let $\epsilon>0$ be sufficiently small, then we set
\[a=\alpha e^{\pi i\epsilon},\quad b=\beta e^{2\pi i\epsilon},\quad
c=\gamma e^{3\pi i\epsilon},\quad d=\delta e^{-6\pi i\epsilon},
\] 
while we keep $t$ undisturbed.
Observe that $|a|=\alpha$, $|b|=\beta$ etc., and that 
$abcd=\alpha\beta\gamma\delta$.
Note that the weight function $W(x)$ (see \eqref{Weight})
corresponding to the parameters $(a,b,c,d,t)$ only has simple poles. 
We can now define the (in general complex) measure $\nu$
with respect to the parameters $(a,b,c,d,t)$ by the same formulas
\eqref{S} and \eqref{nu} as before. We also keep the same notations
for the other parameter-dependent objects we have encountered 
(such as e.g. $\phi_{\gamma}$, $\chi_k$ etc.).
Then for spectral parameters
$\gamma, \gamma^\prime\in {\mathbb{C}}\setminus \{0\}$
with $\gamma^{\pm 1}, (\gamma^\prime)^{\pm 1}\not=q^{1+l}/\tilde{d}$
($l\in {\mathbb{Z}}_+$) we can write for $k\ll 0$,  
\begin{equation}\label{altform}
\int \chi_k(x)\phi_{\gamma}(x)\phi_{\gamma^\prime}(x)d\nu(x)=
\frac{K}{4\pi i}\int_{x\in C_k}\phi_{\gamma}(x)\phi_{\gamma^\prime}(x)\left(
\frac{W(x)}{x}\right)dx,
\end{equation}
where $C_k$ is a continuous, rectifiable Jordan curve in the complex plane, 
satisfying the following additional conditions:
\begin{enumerate}
\item[{\bf --}]
$C_k$ has a parametrization of the form $r_k(z)e^{2\pi iz}$ for 
$z\in [0,1]$ with positive radial function $r_k:[0,1]\rightarrow (0,\infty)$
(we orientate $C_k$ according to this parametrization);
\item[{\bf --}] $C_k$ is invariant under inversion, i.e.
$C_k^{-1}=C_k$, where $C_k^{-1}=\{x^{-1} \,\, | \,\, x\in C_k\}$;
\item[{\bf --}] The sequences 
$\{aq^l, bq^l, cq^l, q^{1+l}/d\}_{l\in {\mathbb{Z}}_+}$ 
are in the interior of $C_k$;
\item[{\bf --}] The intersection of the $q$-interval $I$ with the interior
of $C_k$ is given by the sequence 
$\{dtq^m \,\, | \,\, m\in {\mathbb{Z}}: m\geq k\}$.
\end{enumerate}
The existence of $C_k$ is easy to prove with our choice of perturbed
parameters, and \eqref{altform} is then a direct consequence of 
Cauchy's theorem.

Now recall that the Askey-Wilson function $\phi_{\gamma}$ is an eigenfunction
of the second-order $q$-difference operator $L$ with eigenvalue $\mu(\gamma)$,
see section 3. This implies that
\begin{equation}\label{onehand}
\begin{split}
\bigl(\mu(\gamma)-\mu(\gamma^\prime)\bigr)&\int\chi_k(x)\phi_{\gamma}(x)
\phi_{\gamma^\prime}(x)d\nu(x)=\\
=&\frac{K}{4\pi i}\int_{x\in C_k}\left(
(L\phi_{\gamma})(x)\phi_{\gamma^\prime}(x)-
\phi_{\gamma}(x)\bigl(L\phi_{\gamma^\prime}\bigr)(x)\right)
\left(\frac{W(x)}{x}\right)dx.
\end{split}
\end{equation}
In view of the explicit expression for $L$ 
(see \eqref{LAW}), we can write
\begin{equation}\label{a1}
(L\phi_{\gamma})(x)\phi_{\gamma^\prime}(x)-
\phi_{\gamma}(x)\bigl(L\phi_{\gamma^\prime}\bigr)(x)=
\alpha(1/x)\xi_{\gamma,\gamma^\prime}(x/q)-
\alpha(x)\xi_{\gamma,\gamma^\prime}(x)
\end{equation}
with 
\[
\xi_{\gamma,\gamma^\prime}(x)=\phi_{\gamma}(x)\phi_{\gamma^\prime}(qx)-
\phi_{\gamma}(qx)\phi_{\gamma^\prime}(x). 
\]
Furthermore, by the explicit
expression for the weight function $W(x)$ (see \eqref{Weight}) and
for the function $\alpha(x)$ (see \eqref{LAW}), we have the identity 
\begin{equation}\label{a2}
\alpha(1/x)W(x)=\alpha(x/q)W(x/q).
\end{equation}
Substitution of \eqref{a1} and \eqref{a2} into the right hand side of 
\eqref{onehand} then shows that
\begin{equation}\label{otherhand}
\begin{split}
\bigl(\mu(\gamma)-\mu(\gamma^\prime)\bigr)&\int\chi_k(x)\phi_{\gamma}(x)
\phi_{\gamma^\prime}(x)d\nu(x)=\\
&=\frac{K}{4\pi i}\int_{x\in \,q^{-1}C_k-C_k}
\alpha(x)\xi_{\gamma,\gamma^\prime}(x)\left(\frac{W(x)}{x}\right)dx.
\end{split}
\end{equation}
Using the explicit form of the contour $C_k$, 
we can shrink $q^{-1}C_k$ back to $C_k$ in the right hand side of 
\eqref{otherhand} while picking up residues of the integrand at 
the simple poles $x=dtq^{k-1}$ and $x=d^{-1}t^{-1}q^{-k}$.
By \eqref{a2} and by the identity $\xi_{\gamma,\gamma^\prime}(x)=-
\xi_{\gamma,\gamma^\prime}(q^{-1}x^{-1})$, we can pull the two residues
together, yielding
\begin{equation*}
\begin{split}
\bigl(\mu(\gamma)-\mu(\gamma^\prime)\bigr)&\int\chi_k(x)\phi_{\gamma}(x)
\phi_{\gamma^\prime}(x)d\nu(x)=\\
=&K\underset{y=dtq^{k-1}}{\hbox{Res}}\left(\frac{W(y)}{y}\right)
\alpha(dtq^{k-1})\xi_{\gamma,\gamma^\prime}(dtq^{k-1})\\
=&
\lbrack \phi_{\gamma}, \phi_{\gamma^\prime}\rbrack(dtq^{k-1}),
\end{split}
\end{equation*}
as desired.
\end{proof}


\section{Orthogonality relations.}\label{section7}

If we take in \eqref{kernel} the values of the spectral parameters 
$\gamma, \gamma^\prime$ in the 
discrete support $\widetilde{S}$ of the dual measure $\widetilde{\nu}$, 
then the limit of \eqref{kernel} as $k\rightarrow -\infty$
can be taken point-wise. We evaluate these limits in this section,
which lead to explicit orthogonality relations for the Askey-Wilson functions
$\phi_{\gamma}$ ($\gamma\in\widetilde{S}$) with respect to
the measure $\nu$.
This essentially establishes
theorem \ref{main} for the ``completely discrete part'' of the 
Askey-Wilson function transform.
\begin{lem}\label{proportho}
The Askey-Wilson functions $\phi_{\gamma}$ 
\textup{(}$\gamma\in \widetilde{S}$\textup{)}
are elements of the Hilbert space ${\mathcal{H}}$. Furthermore,
\[
\int\phi_{\gamma}(x)\overline{\phi_{\gamma^\prime}(x)}d\nu(x)=
0,\qquad \gamma,\gamma^\prime\in\widetilde{S},\,\,\, \gamma\not=\gamma^\prime.
\]
\end{lem}
\begin{proof}
We first observe that the asymptotic behaviour of the
discrete weights $\nu\bigl(\{x\}\bigr)$ ($x\in \widetilde{S}_-$) is given
by
\begin{equation}\label{asymptweight}
\nu\bigl(\{dtq^{-k}\}\bigr)=\frac{M}{2\tilde{a}^{2k}}\bigl(1+
\mathcal{O}(q^k)\bigr),\qquad k\rightarrow \infty,
\end{equation}
where $M$ is the positive constant
\begin{equation}\label{M}
M=\frac{\theta(qt)K}{\bigl(q,q;q\bigr)_{\infty}\theta(adt)\theta(bdt)
\theta(cdt)}.
\end{equation}
This is a direct consequence of the explicit expression \eqref{weighttoinfty}
for $\nu\bigl(\{x\}\bigr)$ ($x\in S_-$).

For the moment we assume $ad\not\in q^{\mathbb{Z}}$ and  
$\tilde{a}, \tilde{d}\tilde{t}\not\in\pm q^{\frac{1}{2}{\mathbb{Z}}}$
in order to be able to apply the $c$-function expansion for 
$\phi_{\gamma}$ ($\gamma\in\widetilde{S}$).
We fix $\gamma\in\widetilde{S}$,  so that
$\phi_{\gamma}(x)=\widetilde{c}(\gamma^{-1})\Phi_{\gamma^{-1}}(x)$
for $x\in S_-$ sufficiently negative by proposition \ref{cprop}.
Observe that the apparent poles of $\widetilde{c}(\gamma^{-1})$
as function of the parameters $(a,b,c,d,t)\in V$ are removable.
Hence we can extend the definition of the $c$-function
$\widetilde{c}(\gamma^{-1})$ to parameter
values $(a,b,c,d,t)\in V$ by continuity and we obtain
\begin{equation}\label{half}
\phi_{\gamma}(x)=\widetilde{c}(\gamma^{-1})\Phi_{\gamma^{-1}}(x),\qquad
\gamma\in \widetilde{S}
\end{equation}
for all $(a,b,c,d,t)\in V$ if $x\in I$ is sufficiently negative.
Since 
\begin{equation}\label{Phiasympt}
\Phi_{\gamma^{-1}}(dtq^{-k})=\bigl(\tilde{a}\gamma^{-1}\bigr)^k
\bigl(1+{\mathcal{O}}(q^k)\bigr),\qquad k\rightarrow \infty
\end{equation}
and $|\gamma^{-1}|<1$ for $\gamma\in\widetilde{S}$, 
it follows from \eqref{asymptweight} that
$\phi_{\gamma}\in {\mathcal{H}}$ for all $\gamma\in\widetilde{S}$.

Fix now $\gamma,\gamma^\prime\in \widetilde{S}$ 
with $\gamma\not=\gamma^\prime$.
Then $\mu(\gamma)\not=\mu(\gamma^\prime)$ by \eqref{eigenvalue}, 
hence by proposition 
\ref{propwronskian} and by the previous paragraph,
\begin{equation}\label{b1}
\int\phi_{\gamma}(x)\overline{\phi_{\gamma^\prime}(x)}d\nu(x)=
\frac{1}{(\mu(\gamma)-\mu(\gamma^\prime))}
\lim_{k\rightarrow \infty}\lbrack \phi_{\gamma}, 
\phi_{\gamma^\prime}\rbrack(dtq^{-k}).
\end{equation}
Here we have used 
the fact that $\phi_{\gamma^\prime}(\cdot)$ is real valued
on $\hbox{supp}(\nu)$ in order
to get rid of the complex conjugate in the integrand.

It remains to prove that the limit in the right hand side of \eqref{b1}
tends to zero. Let $f$ and $g$ be two functions in the Hilbert space 
${\mathcal{H}}$. By the asymptotics 
\eqref{asymptweight}, we have $f(dtq^{-k})=o(\tilde{a}^k)$ as 
$k\rightarrow \infty$, and similarly for $g$.
Furthermore, we have the asymptotics
\begin{equation}\label{alphaasympt}
\alpha(dtq^{-k})=\tilde{a}^2\bigl(1+
{\mathcal{O}}(q^k)\bigr),\qquad k\rightarrow\infty
\end{equation}
for the coefficient $\alpha(\cdot)$ of 
the Askey-Wilson second order $q$-difference 
operator $L$ (see \eqref{LAW}).
Combined with \eqref{asymptweight}, it follows from
the definition \eqref{wronskian} of the Wronskian that
\[
\lim_{k\rightarrow\infty}\lbrack f,g\rbrack(dtq^{-k})=0,\qquad \forall
f,g\in {\mathcal{H}}.
\]
In particular, the limit in the right hand side of \eqref{b1} tends to zero,
hence $\phi_{\gamma}$ 
is orthogonal to $\phi_{\gamma^\prime}$ in ${\mathcal{H}}$.
\end{proof}

The quadratic norms of the Askey-Wilson functions $\phi_{\gamma}$
($\gamma\in\widetilde{S}$) in ${\mathcal{H}}$ can be evaluated as
follows.
\begin{lem}\label{normlem}
For all $\gamma\in\widetilde{S}$,
\[\int|\phi_{\gamma}(x)|^2d\nu(x)=
\frac{1}{2\widetilde{\nu}\bigl(\{\gamma\}\bigr)}.
\]
\end{lem}
\begin{proof}
As in the proof of lemma \ref{proportho}, we assume for the moment 
that $ad\not\in q^{\mathbb{Z}}$ and that 
$\tilde{a}, \tilde{d}\tilde{t}\not\in\pm q^{\frac{1}{2}{\mathbb{Z}}}$.
In view of \eqref{half}, \eqref{Phiasympt} and \eqref{asymptweight}, 
these generic conditions can be removed at the end of the proof
by applying the dominated convergence theorem.

We fix $\gamma\in\widetilde{S}$, then by lemma \ref{proportho}
and proposition \ref{propwronskian},
\begin{equation}\label{c1}
\begin{split}
\int |\phi_{\gamma}(x)|^2d\nu(x)=&\lim_{k\rightarrow -\infty}
\int\chi_k(x)\phi_{\gamma}(x)^2d\nu(x)\\
=&\lim_{k\rightarrow -\infty}\left(\lim_{\gamma^\prime\rightarrow \gamma}
\int\chi_k(x)\phi_{\gamma}(x)\phi_{\gamma^\prime}(x)d\nu(x)\right)\\
=&\lim_{k\rightarrow \infty}\left(\lim_{\gamma^\prime\rightarrow\gamma}
\frac{\lbrack \phi_{\gamma}, \phi_{\gamma^\prime}\rbrack(dtq^{-k})}
{(\mu(\gamma)-\mu(\gamma^\prime))}\right),
\end{split}
\end{equation}
where we used that $\phi_{\gamma}(\cdot)$ 
is real valued on $\hbox{supp}(\nu)$ for the first equality.

It remains to evaluate 
the limits of the Wronskian in the last equality of \eqref{c1}.
We use the $c$-function expansion for the Askey-Wilson functions
(see proposition \ref{cprop} and \eqref{half}) to rewrite the
Wronskian as 
\begin{equation}\label{d1}
\begin{split}
\lbrack \phi_{\gamma}, \phi_{\gamma^\prime}\rbrack(dtq^{-k})=
&\widetilde{c}(\gamma^{-1})\widetilde{c}(\gamma^\prime)
\lbrack \Phi_{\gamma^{-1}},\Phi_{\gamma^\prime}\rbrack(dtq^{-k})\\
&+\widetilde{c}(\gamma^{-1})\widetilde{c}(\gamma^\prime{}^{-1})
\lbrack \Phi_{\gamma^{-1}},\Phi_{\gamma^\prime{}^{-1}}\rbrack(dtq^{-k})
\end{split}
\end{equation}
for $k\gg 0$. 
Now there exists open neighbourhoods $U_{\pm}$ of $\gamma^{\pm 1}$
in the complex plane such that
\begin{equation}\label{replacement}
\Phi_{\gamma^\prime}(x)=\Phi_{\gamma^\prime}^{free}(x)(1+f_{\gamma^\prime}(x)),
\qquad x\in S_-,\,\,\, x\ll 0
\end{equation}
where $f_{\gamma^\prime}(x)$ admits a convergent power-series expansion
around $x=-\infty$ with coefficients depending analytically on 
$\gamma^\prime\in U_{\pm}$ and with constant coefficient
equal to zero. Furthermore, for sufficiently small 
neighbourhoods $U_{\pm}$, differentiation  
with respect to $\gamma^\prime$ may be interchanged
with summation in the power series expansion of $f_{\gamma^\prime}(x)$ 
around $x=-\infty$ when $|x|\geq N>0$ for some $U_{\pm}$ independent 
positive constant $N$.

Combined with \eqref{replacement}, 
we see that
\begin{equation}\label{dhalf}
\begin{split}
\lbrack \Phi_{\gamma^{-1}}, &\Phi_{\gamma^\prime{}^{\pm 1}}
\rbrack(dtq^{-k})=\\
=&2\,\nu\bigl(\{dtq^{-k}\}\bigr)\alpha(dtq^{-k})\tilde{a}^{2k-1}
\bigl(\gamma^{-1}\gamma^\prime{}^{\pm 1}\bigr)^{k-1}
\bigl(\gamma^{-1}-\gamma^\prime{}^{\pm 1}\bigr)
\bigl(1+\mathcal{O}(q^k)\bigr)
\end{split}
\end{equation}
as $k\rightarrow \infty$, with $\mathcal{O}(q^k)$ uniform in 
$\gamma^\prime\in U_{\pm}$. 
Now we substitute \eqref{replacement} in \eqref{d1}, and we use
the fact that $\widetilde{c}(\gamma^\prime)$
has a simple zero at $\gamma^\prime=\gamma$ and that
\begin{equation}\label{muterm}
\mu(\gamma)-\mu(\gamma^\prime)=\frac{\tilde{a}}{\gamma}
\bigl(\gamma-\gamma^\prime\bigr)\bigl(\gamma-\gamma^\prime{}^{-1}\bigr),
\end{equation}
then using \eqref{dhalf} and the fact that $|\gamma|>1$, we derive that
\begin{equation}\label{e1}
\lbrack \phi_{\gamma}, \phi_{\gamma^\prime}\rbrack(dtq^{-k})=
2\nu\bigl(\{dtq^{-k}\}\bigr)\alpha(dtq^{-k})
\tilde{a}^{2k-1}\bigl(s_k(\gamma,\gamma^\prime)+r_k(\gamma,\gamma^\prime)\bigr)
\end{equation}
with
\begin{equation}\label{sterm}
s_k(\gamma,\gamma^\prime)=\sum_{\epsilon=\pm 1}
\widetilde{c}(\gamma^\prime{}^\epsilon)\widetilde{c}(\gamma^{-1})
\gamma^{\prime}{}^{\epsilon(k-1)}\gamma^{1-k}
\bigl(\gamma^{-1}-\gamma^\prime{}^\epsilon\bigr)
\end{equation}
and with a remainder term $r_k(\gamma,\gamma^\prime)$ satisfying 
\begin{equation}\label{remainder}
\lim_{\gamma^\prime\rightarrow \gamma}\frac{r_k(\gamma,\gamma^\prime)}
{(\mu(\gamma)-\mu(\gamma^\prime))}={\mathcal{O}}(q^k),\qquad k\rightarrow 
\infty.
\end{equation}
Using \eqref{muterm}, \eqref{sterm} and the fact that 
$\widetilde{c}(\gamma^\prime)$ has a simple zero at $\gamma^\prime=\gamma$,
we derive
\begin{equation}\label{e2}
\begin{split}
\lim_{\gamma^\prime\rightarrow\gamma}\frac{s_k(\gamma,\gamma^\prime)}
{(\mu(\gamma)-\mu(\gamma^\prime))}
=&\frac{\gamma}{\tilde{a}}\left(\frac{1}{\widetilde{c}(\gamma^{-1})}
\underset{\gamma^\prime=\gamma}{\hbox{Res}}\left(
\frac{1}{\widetilde{c}(\gamma^\prime)}\right)\right)^{-1}+
\frac{\gamma}{\tilde{a}}\frac{\widetilde{c}(\gamma^{-1})^2}
{(\gamma^{-1}-\gamma)\gamma^{2k}}\\
=&\frac{\widetilde{K}}{\tilde{a}\widetilde{c}_0}
\frac{1}{2\widetilde{\nu}\bigl(\{\gamma\}\bigr)}
+\frac{\gamma}{\tilde{a}}\frac{\widetilde{c}(\gamma^{-1})^2}
{(\gamma^{-1}-\gamma)\gamma^{2k}}
\end{split}
\end{equation}
where $\widetilde{K}$ and $\widetilde{c}_0$ are the constants $K$
(see \eqref{K}) and $c_0$ (see \eqref{c0funct}) respectively with respect to
dual parameters. Here we have used \eqref{Weight} and \eqref{nu} 
for the second equality. 

We can evaluate now the
limits in the last equality of \eqref{c1} 
by substituting \eqref{e1} for the Wronskian and by
using \eqref{asymptweight}, \eqref{alphaasympt}, 
\eqref{sterm} and \eqref{remainder} together with the fact that $|\gamma|>1$.
This gives
\[\int|\phi_{\gamma}(x)|^2d\nu(x)=
\frac{\widetilde{K}M}{\widetilde{c}_0}\frac{1}
{2\widetilde{\nu}\bigl(\{\gamma\}\bigr)}.
\]
By a direct computation one verifies that
$\widetilde{K}M/\widetilde{c}_0=1$, which completes the proof of the lemma.
\end{proof}

Let ${\mathcal{D}}_d\subset {\mathcal{D}}$ be the sub-space consisting
of the functions in ${\mathcal{D}}$ which are supported on the discrete 
support $S$ of the measure $\nu$, and let $\widetilde{\mathcal{D}}_d$
be
the sub-space ${\mathcal{D}}_d$ with respect to dual parameters.
In view of lemma \ref{rewrite} and lemma \ref{dualdomain},
we can now reformulate lemma \ref{proportho}
and lemma \ref{normlem} in the following way.

\begin{prop}\label{dd}
Let $f,g\in {\mathcal{D}}_d$, then
${\mathcal{F}}f, {\mathcal{F}}g\in \widetilde{\mathcal{H}}$ and 
\[\int\bigl({\mathcal{F}}f\bigr)(\gamma)
\overline{\bigl({\mathcal{F}}g\bigr)(\gamma)}d\widetilde{\nu}(\gamma)=
\int f(x){\overline{g(x)}}d\nu(x).
\]
\end{prop}

Recall from \eqref{redpol} that $\phi_{\gamma}$ is a multiple
of the Askey-Wilson polynomial $p_k$ when 
$\gamma=\gamma_k=\tilde{a}q^k\in \widetilde{S}_+$ ($k\in {\mathbb{Z}}_+$). 
So if $\tilde{a}>1$, then
lemma \ref{proportho} and lemma \ref{normlem} give explicit orthogonality
relations for the ``low degree'' Askey-Wilson polynomials
$p_k$ ($k\in {\mathbb{Z}}_+$, $k\leq k_0$) with respect to
the one-parameter family of measures $\nu$ ($t<0$), where $k_0$
is the largest positive integer such that $\tilde{a}q^{k_0}>1$.
The sub-space spanned by $\{p_k\}_{k=0}^{k_0}$
is exactly the sub-space of polynomials in $x+x^{-1}$ which are square
integrable with respect to $\nu$ (this follows easily from
\eqref{asymptweight}).
Furthermore, the functions $\phi_{\gamma}\in {\mathcal{H}}$ with 
$\gamma\in \widetilde{S}_-$ 
constitute an explicit family of mutually orthogonal
functions with respect to the measure $\nu$, which in addition
are orthogonal to the finite set of Askey-Wilson polynomials $p_k$
($k\in {\mathbb{Z}}_+$, $k\leq k_0$).

{}From this point of view, lemma \ref{proportho}
and lemma \ref{normlem} bear close resemblance with 
indeterminate moment problems and non-extremal measures. 
In fact, 
the measures $\nu$ ($t<0$) formally reduce to a 
genuine one-parameter family
of non-extremal orthogonality measures for some indeterminate
moment problem in certain degenerate cases 
of the Askey-Wilson function transform. 
This point of view was emphasized 
in \cite{KS3}, where it was for instance shown that a formal limit
of the above orthogonality relations leads to a one-parameter family of 
non-extremal orthogonality 
measures for the continuous dual $q^{-1}$-Hahn polynomials,
see also \cite{KS1}.

\begin{rem}\label{Bessel}
Fix $k\in {\mathbb{Z}}$ sufficiently negative, and let 
$\gamma_1(k), \gamma_2(k)$ be two zeros of the function 
$\gamma\mapsto \phi_{\gamma}(dtq^{k-1})$ such that 
$\mu(\gamma_1(k))\not=\mu(\gamma_2(k))$. 
Then proposition \ref{propwronskian} implies that the Askey-Wilson
functions $\phi_{\gamma_1(k)}$
and $\phi_{\gamma_2(k)}$ are mutually orthogonal with respect to
the  measure $\chi_k(x)d\nu(x)$ of {\it compact} support. 
These Fourier-Bessel type
orthogonality relations for the Askey-Wilson functions 
were derived and studied by Suslov in \cite{S} and \cite{S2}.
\end{rem}


\section{The continuous part of the Askey-Wilson function transform.}
\label{section8} 

In this section we first evaluate the weak limit of \eqref{kernel} as 
$k\rightarrow -\infty$ with respect to the product
measure $\widetilde{\nu}|_{\mathbb{T}}\times \widetilde{\nu}|_{\mathbb{T}}$.
We directly formulate the result 
in terms of suitable unitarity properties of the Askey-Wilson 
function transform (compare with proposition \ref{dd} in the square
integrable setting).

Let ${\mathcal{D}}_c\subset {\mathcal{D}}$ be the sub-space consisting
of the functions in ${\mathcal{D}}$ which are
supported on ${\mathbb{T}}$, and let $\widetilde{\mathcal{D}}_c$
be the sub-space ${\mathcal{D}}_c$ with respect to dual parameters.
Observe that ${\mathcal{D}}=
{\mathcal{D}}_c\oplus {\mathcal{D}}_d$ 
is an orthogonal direct sum decomposition of the dense
sub-space ${\mathcal{D}}\subset {\mathcal{H}}$. 

\begin{prop}\label{cc}
Let $f,g\in {\mathcal{D}}_c$, 
then ${\mathcal{F}}f, {\mathcal{F}}g \in \widetilde{\mathcal{H}} $
and 
\[\int \bigl({\mathcal{F}}f\bigr)(\gamma)
\overline{\bigl({\mathcal{F}}g\bigr)(\gamma)}d\widetilde{\nu}(\gamma)=
\int f(x)\overline{g(x)}d\nu(x).
\]
\end{prop}
\begin{proof}
The proof is similar to the proof of \cite[prop 7.7]{KS1}
and of \cite[prop 6.1]{KMU}, where the analogous statement
was derived for the big and the little $q$-Jacobi function transforms,
respectively. Since some care has to be taken in order
to match the constants, we repeat here the proof in some detail.
 
We prove the proposition with respect to dual parameters 
(cf. lemma \ref{dualdomain}). It suffices to prove it for functions
$f,g\in \widetilde{\mathcal{D}}_c$ which are continuous on $\mathbb{T}$
and supported within $\mathbb{T}\setminus \{\pm 1\}$. 

We use lemma \ref{rewrite}, proposition \ref{propwronskian} and the
explicit form of the measure $\widetilde{\nu}|_{\mathbb{T}}$ (see 
\eqref{Weight} and \eqref{nu}), together
with the invariance of the Askey-Wilson function $\phi_{\gamma}$
and the measure $d\widetilde{\nu}(\gamma)$ under $\gamma\mapsto
\gamma^{-1}$, to write
\begin{equation}\label{g1}
\begin{split}
\int &\chi_{1-k}(x)\bigl(\widetilde{\mathcal{F}}f\bigr)(x)
\overline{\bigl(\widetilde{\mathcal{F}}g\bigr)(x)}d\nu(x)=\\
&=\frac{\widetilde{K}^2}{4\pi^2\widetilde{c}_0^2}
\int_{\theta^\prime=0}^{\pi}\int_{\theta=0}^{\pi}
f(e^{i\theta})\overline{g(e^{i\theta^\prime})}
\frac{\lbrack \phi_{e^{i\theta}}, \phi_{e^{i\theta^\prime}}\rbrack(dtq^{-k})}
{\bigl(\mu(e^{i\theta})-\mu(e^{i\theta^\prime})\bigr)}
\frac{d\theta}{|\widetilde{c}(e^{i\theta})|^2} 
\frac{d\theta^\prime}{|\widetilde{c}(e^{i\theta^\prime})|^2}
\end{split}
\end{equation}
for $k\in {\mathbb{Z}}$ sufficiently large. It suffices to prove that
the right hand side of \eqref{g1}
tends to $\int f(\gamma)\overline{g(\gamma)}d\widetilde{\nu}(\gamma)$
as $k\rightarrow \infty$.

In order to compute the limit $k\rightarrow \infty$ of the right hand
side of \eqref{g1}, we first observe that
the factor $\mu(e^{i\theta})-\mu(e^{i\theta^\prime})$ occuring in
the integrand can be rewritten as
\begin{equation}\label{g2}
\begin{split}
\mu(e^{i\theta})-\mu(e^{i\theta^\prime})&=
\tilde{a}e^{-i\theta}\bigl(e^{i\theta}-e^{i\theta^\prime}\bigr)
\bigl(e^{i\theta}-e^{-i\theta^\prime}\bigr)\\
&=-4\tilde{a}\sin\left(\frac{\theta+\theta^\prime}{2}\right)
\sin\left(\frac{\theta-\theta^\prime}{2}\right),
\end{split}
\end{equation}
cf. \eqref{muterm}. Similarly as in 
the proof of lemma \ref{normlem}, we derive now from 
the first equality in \eqref{g2}, the $c$-function
expansion (see proposition \ref{cprop}) and the
mean value theorem, that
\begin{equation}\label{hulpnew}
 \lbrack \phi_{e^{i\theta}}, \phi_{e^{i\theta^\prime}}\rbrack(dtq^{-k})
=2\nu\bigl(\{dtq^{-k}\}\bigr)\alpha(dtq^{-k})\tilde{a}^{2k-1}
\bigl(s_k(\theta,\theta^\prime)+r_k(\theta,\theta^\prime)\bigr)
\end{equation}
with 
\begin{equation}\label{sk}
s_k(\theta,\theta^\prime)=
\sum_{\epsilon, \epsilon^\prime=\pm 1}
\widetilde{c}(e^{i\epsilon\theta})\widetilde{c}(e^{i\epsilon^\prime
\theta^\prime})e^{i(k-1)\epsilon\theta}e^{i(k-1)\epsilon^\prime\theta^\prime}
\bigl(e^{i\epsilon\theta}-e^{i\epsilon^\prime\theta^\prime}\bigr)
\end{equation}
and with a remainder term $r_k(\theta,\theta^\prime)$ satisfying
\begin{equation}\label{hulpnew2}
\underset{\delta\leq\theta\not=\theta^\prime\leq \pi-\delta}{\hbox{sup}}
\left|\frac{r_k(\theta,\theta^\prime)}{\bigl(\mu(e^{i\theta})-
\mu(e^{i\theta^\prime})\bigr)}\right|={\mathcal{O}}(kq^k),
\qquad k\rightarrow \infty
\end{equation}
for all $0<\delta<\pi/2$.
We substitute the expression \eqref{hulpnew} for the Wronskian
in the right hand side of \eqref{g1}. In view of the asymptotics
\eqref{asymptweight} and \eqref{alphaasympt}, it then 
suffices to calculate the limit $k\rightarrow \infty$ of
\begin{equation}\label{ghalf}
 \frac{\widetilde{K}^2M\tilde{a}}{4\pi^2\widetilde{c}_0^2}
\int_{\theta^\prime=0}^{\pi}\int_{\theta=0}^\pi
f(e^{i\theta})\overline{g(e^{i\theta^\prime})}
\frac{\bigl(s_k(\theta,\theta^\prime)+r_k(\theta, \theta^\prime)\bigr)}
{\bigl(\mu(e^{i\theta})-\mu(e^{i\theta^\prime})\bigr)}
\frac{d\theta}{|\widetilde{c}(e^{i\theta})|^2}
\frac{d\theta^\prime}{|\widetilde{c}(e^{i\theta^\prime})|^2}.
\end{equation}
The asymptotic behaviour \eqref{hulpnew2} 
shows that the limit $k\to\infty$ of the integral \eqref{ghalf}
of the remainder term $r_k(\theta,\theta^\prime)$
gives zero by dominated convergence. 
By the Riemann-Lebesgue lemma,  the integral \eqref{ghalf} of the two 
terms from $s_k(\theta,\theta^\prime)$ corresponding
to $\epsilon\epsilon^\prime=1$ in \eqref{sk}
tend to zero as $k\to\infty$. 
It thus remains to calculate the limit $k\rightarrow\infty$ of
\eqref{ghalf} in which the factor $\bigl(s_k(\theta,\theta^\prime)+
r_k(\theta,\theta^\prime)\bigr)$ of the integrand is replaced by
\begin{equation*}
\begin{split}
t_k(\theta,\theta^\prime)&=\sum_{\stackrel{
\scriptstyle{\epsilon, \epsilon^\prime=\pm 1}}
{\epsilon\epsilon^\prime=-1}}
\widetilde{c}(e^{i\epsilon\theta})\widetilde{c}(e^{i\epsilon^\prime
\theta^\prime})e^{i(k-1)\epsilon\theta}e^{i(k-1)\epsilon^\prime\theta^\prime}
\bigl(e^{i\epsilon\theta}-e^{i\epsilon^\prime\theta^\prime}\bigr)\\
&=\bigl(\widetilde{c}(e^{-i\theta})\widetilde{c}(e^{i\theta^\prime})
-\widetilde{c}(e^{i\theta})\widetilde{c}(e^{-i\theta^\prime})\bigr)
e^{i(k-1)(\theta-\theta^\prime)}\bigl(e^{-i\theta}-e^{i\theta^\prime}\bigr)\\
&-4\sin\left(\frac{\theta+\theta^\prime}{2}\right)\sin\left(
\frac{(2k-1)(\theta-\theta^\prime)}{2}\right)\widetilde{c}(e^{i\theta})
\widetilde{c}(e^{-i\theta^\prime}).
\end{split}
\end{equation*}
It follows now by yet another application of the Riemann-Lebesgue lemma 
that it remains to calculate the limit 
$k\rightarrow\infty$ of 
\begin{equation}\label{g4}
\frac{\widetilde{K}^2M}{4\pi^2\widetilde{c}_0^2}
\int_{\theta^\prime=0}^\pi\int_{\theta=0}^\pi
f(e^{i\theta})\overline{g(e^{i\theta^\prime})}
\left(\frac{\sin\left((2k-1)(\theta-\theta^\prime)/2\right)}
{\sin\left((\theta-\theta^\prime)/2\right)}\right)\frac{d\theta}
{\widetilde{c}\bigl(e^{-i\theta}\bigr)}
\frac{d\theta^\prime}{\widetilde{c}\bigl(e^{i\theta^\prime}\bigr)}.
\end{equation}
Now by the well known $L^2$-properties of the Dirichlet kernel, the
limit $k\rightarrow \infty$ of \eqref{g4} exists,
and it equals
\[\frac{\widetilde{K}^2M}{2\pi\widetilde{c}_0^2}
\int_{\theta=0}^\pi f(e^{i\theta})\overline{g(e^{i\theta})}
\frac{d\theta}{|\widetilde{c}(e^{i\theta})|^2}
=\frac{\widetilde{K}M}{\widetilde{c}_0}
\int f(\gamma)\overline{g(\gamma)}d\widetilde{\nu}(\gamma).
\]
The proposition follows now directly from the fact that
$\widetilde{K}M/\widetilde{c}_0=1$ (compare with the proof of
lemma \ref{normlem}).
\end{proof}

In order to completely understand the weak limit of \eqref{kernel}
as $k\rightarrow -\infty$ with respect to the product measure
$\widetilde{\nu}\times \widetilde{\nu}$, we still have to deal
with the mixed continuous-discrete case. 
This case is covered by the following lemma.

\begin{lem}\label{dc}
Let $f\in {\mathcal{D}}_c$ and $g\in {\mathcal{D}}_d$,
then
\[\int \bigl({\mathcal{F}}f\bigr)(\gamma)
\overline{\bigl({\mathcal{F}}g\bigr)(\gamma)}
d\widetilde{\nu}(\gamma)=0.
\]
\end{lem}
\begin{proof}
We establish the desired identity with respect to dual parameters.
Let $f\in \widetilde{\mathcal{D}}_c$ and $g\in \widetilde{\mathcal{D}}_d$.
By proposition \ref{dd} and proposition \ref{cc} we have 
$\widetilde{\mathcal{F}}f\in \mathcal{H}$ and 
$\widetilde{\mathcal{F}}g\in \mathcal{H}$ respectively.
In particular, we may assume without loss of generality that
$f$ is continuous on ${\mathbb{T}}$, and supported within 
${\mathbb{T}}\setminus \{\pm 1\}$. For $\gamma\in \widetilde{S}\cup
\widetilde{S}^{-1}$ we now define $f_{\gamma}\in \widetilde{\mathcal{D}}_c$ by
$f_{\gamma}(\gamma^\prime)=f(\gamma^\prime)/\bigl(\mu(\gamma)-
\mu(\gamma^\prime)\bigr)$. Observe that $f_{\gamma}$ is continuous on 
${\mathbb{T}}$ and supported within ${\mathbb{T}}\setminus \{\pm 1\}$.
By lemma \ref{rewrite} and proposition \ref{propwronskian} we can now write
\begin{equation}\label{f1}
\begin{split}
\int\bigl(\widetilde{\mathcal{F}}f\bigr)(x)
&\overline{\bigl(\widetilde{\mathcal{F}}g\bigr)(x)}d\nu(x)=\\
&=2\sum_{\gamma\in\widetilde{S}}
\overline{g(\gamma)}\widetilde{\nu}(\{\gamma\})
\lim_{k\rightarrow \infty}
\int_{\mathbb{T}} f_{\gamma}(\gamma^\prime)\lbrack \phi_{\gamma}, 
\phi_{\gamma^\prime}\rbrack(dtq^{-k})d\widetilde{\nu}(\gamma^\prime).
\end{split}
\end{equation}

Now we observe that \eqref{half}, \eqref{Phiasympt}, 
proposition \ref{cprop} and 
the asymptotic behaviour 
\[\underset{\gamma^\prime\in supp(f_{\gamma})}
{\hbox{sup}}\left|\Phi_{\gamma^\prime{}^{\pm 1}}(dtq^{-k})\right|=
{\mathcal{O}}\bigl(\tilde{a}^k\bigr),\qquad k\rightarrow \infty
\]
imply that
\begin{equation*}
\begin{split}
\underset{\gamma^\prime\in supp(f_{\gamma})}
{\hbox{sup}}\left|\lbrack \phi_{\gamma}, 
\phi_{\gamma^\prime}\rbrack(dtq^{-k})\right|&=
2\nu\bigl(\{dtq^{-k}\}\bigr)\alpha(dtq^{-k})\tilde{a}^{2k-1}|\gamma|^{-k}
{\mathcal{O}}(1)\\
&={\mathcal{O}}\bigl(|\gamma|^{-k}\bigr)
\end{split}
\end{equation*}
as $k\rightarrow \infty$, where the second equality is a consequence of
\eqref{asymptweight} and \eqref{alphaasympt}.
But $|\gamma|^{-1}<1$ for $\gamma\in\widetilde{S}$, hence we may
use Lebesgue's dominated convergence theorem to interchange 
limit and integration in the right hand side of \eqref{f1}. It follows
that the right hand side of 
\eqref{f1} tends to zero, which completes the proof of the lemma.
\end{proof}

\section{Completion of the proof of theorem \ref{main}.}

In this section we complete the proof of theorem 
\ref{main}. 

First of all, we observe that the results of section \ref{section7}
and section \ref{section8} immediately imply that the Askey-Wilson 
function transform ${\mathcal{F}}$ extends to an isometry
${\mathcal{F}}: {\mathcal{H}}\rightarrow \widetilde{\mathcal{H}}$.
Indeed, it follows from proposition \ref{dd}, proposition \ref{cc}
and lemma \ref{dc}, together with the orthogonal direct sum decomposition
${\mathcal{D}}={\mathcal{D}}_c\oplus {\mathcal{D}}_d$, that
\[
\int \bigl({\mathcal{F}}f\bigr)(\gamma)
\overline{\bigl({\mathcal{F}}g\bigr)(\gamma)}d\widetilde{\nu}(\gamma)=
\int f(x)\overline{g(x)}d\nu(x),\qquad f,g\in {\mathcal{D}}.
\]
Since ${\mathcal{D}}\subset {\mathcal{H}}$ is dense, it follows that
the Askey-Wilson function transform 
${\mathcal{F}}$ uniquely extends to an isometry 
${\mathcal{F}}: {\mathcal{H}}\rightarrow \widetilde{\mathcal{H}}$
by continuity. In particular, the dual Askey-Wilson function transform 
$\widetilde{\mathcal{F}}: \widetilde{\mathcal{H}}\rightarrow {\mathcal{H}}$
is an isometry in view of lemma \ref{dualdomain}.

Fix now arbitrary $f\in \widetilde{\mathcal{H}}$ and $g\in {\mathcal{H}}$,
and write $\widetilde{\chi}_k\in \widetilde{\mathcal{D}}$ 
($k\in {\mathbb{Z}}$ sufficiently negative) 
for the characteristic function $\chi_k$ with
respect to dual parameters (see the beginning of section \ref{WWW}). 
Since ${\mathcal{F}}: {\mathcal{H}}\rightarrow
\widetilde{\mathcal{H}}$  and $\widetilde{\mathcal{F}}: \widetilde{\mathcal{H}}
\rightarrow {\mathcal{H}}$ are continuous, we have
\begin{equation*}
\begin{split}
\int\bigl(\widetilde{\mathcal{F}}f\bigr)(x)\overline{g(x)}d\nu(x)&=
\lim_{k,m\rightarrow -\infty}
\int \left(\int 
\widetilde{\chi}_k(\gamma)f(\gamma)\overline{\widetilde{\phi}_x(\gamma)}
d\widetilde{\nu}(\gamma)\right)\overline{\chi_m(x)g(x)}d\nu(x)\\
&=\lim_{k,m\rightarrow -\infty}\int \widetilde{\chi}_k(\gamma)f(\gamma)
\overline{\left(\int \chi_m(x)g(x)\overline{\phi_{\gamma}(x)}d\nu(x)\right)}
d\widetilde{\nu}(\gamma)\\
&=\int f(\gamma)\overline{\bigl({\mathcal{F}}g\bigr)(\gamma)}
d\widetilde{\nu}(\gamma),
\end{split}
\end{equation*}
where we used Fubini's theorem, 
the duality \eqref{dualityfunction} of the Askey-Wilson function,
and the fact that the Askey-Wilson function $\phi_{\gamma}(x)$
is real-valued for $x\in \hbox{supp}(\nu)$ and
$\gamma\in \hbox{supp}(\widetilde{\nu})$. 
But then we have for all $f,g \in {\mathcal{H}}$,
\[
\int \bigl(\widetilde{\mathcal{F}}\bigl({\mathcal{F}}f\bigr)\bigr)(x)
\overline{g(x)}d\nu(x)=
\int \bigl({\mathcal{F}}f\bigr)(\gamma)
\overline{\bigl({\mathcal{F}}g\bigr)(\gamma)}d\widetilde{\nu}(\gamma)
=\int f(x)\overline{g(x)}d\nu(x)
\]
since ${\mathcal{F}}: {\mathcal{H}}\rightarrow \widetilde{\mathcal{H}}$
is an isometry. It follows that
$\widetilde{\mathcal{F}}\circ {\mathcal{F}}=\hbox{Id}_{\mathcal{H}}$.
By a similar argument, we obtain 
${\mathcal{F}}\circ \widetilde{\mathcal{F}}=
\hbox{Id}_{\widetilde{\mathcal{H}}}$ (or simply replace the parameters 
by the corresponding dual parameters
in $\widetilde{\mathcal{F}}\circ {\mathcal{F}}=\hbox{Id}_{\mathcal{H}}$
and use lemma \ref{dualdomain}). 

We conclude that
${\mathcal{F}}: {\mathcal{H}}\rightarrow \widetilde{\mathcal{H}}$
and $\widetilde{\mathcal{F}}: \widetilde{\mathcal{H}}\rightarrow
{\mathcal{H}}$ are isometric isomorphisms, and that $\widetilde{\mathcal{F}}=
{\mathcal{F}}^{-1}$. This completes the proof of theorem \ref{main}.


\bibliographystyle{amsplain}

\end{document}